\documentclass[12pt]{amsart}
\usepackage{fullpage}
\usepackage{hyperref}
\usepackage{amsrefs}
\usepackage{xy}
\xyoption{all}

\newtheorem{thm}{Theorem}[section]
\newtheorem{cor}[thm]{Corollary}
\newtheorem*{thm*}{Theorem}
\newtheorem{prop}[thm]{Proposition}
\theoremstyle{definition} \newtheorem*{dfn*}{Definition}
\theoremstyle{remark} \newtheorem{rem}[thm]{Remark}
\numberwithin{equation}{subsection}

\newcommand{\bbR}{\mathbb{R}} \newcommand{\bbC}{\mathbb{C}} \newcommand{\bbH}{\mathbb{H}}
\newcommand{\bbZ}{\mathbb{Z}} \newcommand{\bbN}{\mathbb{N}}
 \newcommand{\mcO}{\mathcal{O}} 
\newcommand{\mcP}{\mathcal{P}} \newcommand{\mcH}{\mathcal{H}}
\newcommand{\mcS}{\mathcal{S}} \newcommand{\mcU}{\mathcal{U}}

\renewcommand{\H}{\mbox{HS}}
\newcommand{\cG}{{\overline{G}}}
\newcommand{\cK}{{\overline{K}}}
\newcommand{\cmcS}{{\overline{\mcS}}}
\newcommand{\cW}{{\widetilde{W}}}
 \newcommand{\fg}{\mathfrak g} \newcommand{\gl}{\mathfrak{gl}} \renewcommand{\sp}{\mathfrak{sp}} \newcommand{\so}{\mathfrak{so}}
\newcommand{\twoline}[2]{\genfrac{}{}{0pt}{}{#1}{#2}}
\newcommand{\la}{\lambda}
\renewcommand{\a}{\alpha}
\renewcommand{\b}{\beta}
\newcommand{\g}{\gamma}
\newcommand{\vx}{\bar{\text{x}}} \newcommand{\vy}{\bar{\text{y}}} \newcommand{\vz}{\bar{\text{z}}}
\newcommand{\ch}{\text{ch}}
\newcommand{\F}[2]{F^{#1}_{(#2)}} \renewcommand{\c}[3]{{c^{#1}_{{#2}{#3}}}} \newcommand{\ot}{\otimes} \newcommand{\tot}{\widehat \otimes} \newcommand{\dual}[1]{ {\left({#1}\right)}^*}

\begin{document}

\title[Stable Hilbert Series]{Stable Hilbert series of $\mathcal S(\mathfrak g)^K$ for classical groups}
\author{Jeb F. Willenbring}
\email{jw@uwm.edu}
\date{November 3, 2005}
\thanks{The author was supported in part by NSA Grant \# H98230-05-1-0078.}

\address{Jeb F. Willenbring, University of Wisconsin-Milwaukee, Dept. of Math. Sciences, P.O. Box 413, Milwaukee, WI 53201}
\begin{abstract}
Given a classical symmetric pair, $(G,K)$, with $\mathfrak g =
Lie(G)$, we provide descriptions of the Hilbert series of the
algebra of $K$-invariant vectors in the associated graded algebra of
$\mathcal U(\mathfrak g)$ viewed as a $K$-representation under
restriction of the adjoint representation. The description
illuminates a certain stable behavior of the Hilbert series, which
is investigated in a case-by-case basis.  We note that the stable
Hilbert series of one symmetric pair often coincides with others.
Also, for the case of the real form $U(p,q)$ we derive a closed
expression for the Hilbert series when $\min(p,q) \rightarrow
\infty$.
\end{abstract}
\maketitle

\section{Introduction.}$\empty$
The adjoint representation of a reductive group $G$ on its Lie
algebra $\fg$ induces an action on the universal enveloping algebra
$\mcU(\fg)$. The structure of $\mcU(\fg)$ as a representation of $G$
follows from Kostant's theory of harmonic polynomials (see
\cite{Kost63}) on $\fg$. We approach the problem of understanding
the restriction of the induced $Ad(G)$-action on $\mcU(\fg)$ to a
symmetric subgroup $K$. In particular, we focus on the $K$-invariant
subalgebra, $\mcU(\fg)^K$. It is well established that $\mcU(\fg)^K$
has an extraordinarily complicated ring structure, while at the same
time, of fundamental importance to the structure of admissible
representations of real reductive groups.  We present results on the
stability of the Hilbert series of the associated graded
$K$-invariant subalgebra. Namely, we provide results that the stable
Hilbert series of one symmetric pair coincides with that of another.
As a consequence of the PBW theorem, we have $\mcU(\fg) \cong
\mcS(\fg)$ (the associated graded algebra) as a $G$-representation.
Note that here $\mcS(\cdot)$ denotes the symmetric algebra.

In light of this fact we will consider the algebra $\mcS(\fg)^K$,
which even though commutative, also evades a simple presentation in
terms of generators and relations in most examples. Specifically,
consider classical symmetric pairs.  Equivalently, we consider the
situation where $\fg$ is the complexification of the Lie algebra of
the real form, $G_0$, of a classical linear algebraic group over
$\bbC$. The group $K$ will be the complex group corresponding to the
maximal compact subgroup of $G_0$.  The ten families that we
consider are described in Section \ref{subsec_realforms}.

This paper consists of four sections.  In this first section we set
up basic background material and notation in order to state the main
result in Section \ref{subsec_main_thm}.  In Section
\ref{subsec_stability} we interpret the main theorem as a ``stable
limit".  In Section \ref{sec_background} we recall the necessary
results to establish the main theorem, which is proved in Section
\ref{sec_main_thm_proof}.  In Section \ref{sec_closed_expressions}
we show how one can obtain a closed expression for the Hilbert
series of $\mcU(\fg)^K$ in the case of $U(p,q)$ when $\min(p,q)
\rightarrow \infty$.

In order to state the main theorem we set up some standard notation
and review relevant results.  These included:  Symmetric pairs and
real forms (\ref{subsec_realforms}), Hilbert series
(\ref{subsec_hilbert_series}), Partitions
(\ref{subsec_combinatorics}), Irreducible representations of
$GL_n(\bbC)$ (\ref{subsec_irreps}), Multiplicity
(\ref{subsec_multiplicities}), and Littlewood-Richardson
coefficients (\ref{subsec_LRC}).

\subsection{Notation for symmetric pairs and the real forms of the classical
groups.}\label{subsec_realforms} We work in the context of linear
algebraic groups over the complex numbers, but the motivation for
the problem comes from the real forms of these groups. To avoid any
possible ambiguity, we recall the relevant background material,
which can be found in \cite{GW98}, Chapters 1, 11 and 12. To begin
with, a linear algebraic group, $G$, is a Zariski closed subgroup of
$GL_k(\bbC)$.  We let $\mcO(G)$ denote the algebra of regular
functions on $G$.  In this paper, the most general such groups we
consider are (products of) the three classical groups $GL_n(\bbC)$,
$O_n(\bbC)$ and $Sp_{2m}(\bbC)$.  Unless otherwise stated we take
the forms for $O_n(\bbC)$ and $Sp_{2m}(\bbC)$ to be the usual $(x,y)
= x y^t$ and $(x,y) = x J_m y^t$ respectively, where:
\[
    J_m =
    \left[ \begin{array}{cc}   0 & I_m \\ -I_m &   0 \end{array} \right],
\] with $I_m$ being the $m \times m$ identity matrix.

By a \emph{symmetric pair} we main an ordered pair $(G,K)$ such that
$G$ is reductive complex linear algebraic group and $K$ consists of
the fixed points of a \emph{regular} involution, $\theta$, of $G$.
Consequently, the group $K$ will be a complex reductive linear
algebraic group.

If $G$ is a complex algebraic group, a \emph{complex conjugation} on
$G$ is an abstract group involution $\tau$, of $G$ such that $f^\tau
\in \mcO(G)$ for all $f \in \mcO(G)$, where $f^\tau$ is the (complex
valued) function on $G$ defined by $f^\tau(g) =
\overline{f(\tau(g))}$, for all $g \in G$ (see \cite{GW98} Section
1.4.2).

A \emph{real form} of a complex reductive linear algebraic group $G$
is a subgroup $G_0 \subseteq G$ such that there exists a complex
conjugation, $\tau$, such that $G_0 = \{ g \in G | \tau(g) = g
\mbox{ for all } g \in G \}$.  In general, a real form does not have
the structure of a complex linear algebraic group, however it is a
standard fact that the real forms of $G$ correspond to the
symmetrically embedded subgroups of $G$ (see \cite{GW98}).  If $G_0$
is a real form of $G$, we $G$ is the \emph{complexification} of
$G_0$, while $K$ is the complexification of a maximal compact
subgroup, $K_0$.

For our purposes, we will only be interested in the relevant
symmetric pair corresponding to a real form $G_0$.  In order to set
up a convenient notation, we will refer to the symmetric pairs by
``standard" notation for the corresponding real form as in the three
tables below:
\[
\begin{array}{|c||c|c|c|c|} \hline
   G_0     &                  G                  &                      K             &  \theta: G \rightarrow G            \\ \hline \hline
 U(p,q)    &            GL_{p+q}(\bbC)           &    GL_p(\bbC) \times GL_q(\bbC)    &  g \mapsto I_{p,q} g I_{p,q}        \\ \hline
GL(n,\bbR) &              GL_n(\bbC)             &             O_n(\bbC)              &  g \mapsto (g^{-1})^t               \\ \hline
GL(m,\bbH) &            GL_{2m}(\bbC)            &           Sp_{2m}(\bbC)            &  g \mapsto -J_m \, (g^{-1})^t \, J_m \\ \hline
\end{array}
\] where
\[
    I_{p,q} =
    \left[ \begin{array}{cc} I_p & 0   \\    0 & I_q \end{array}
    \right].
\]
In the first column in the above table are three real forms of the
complex general linear group.  We have indicated the corresponding
symmetric pair in the second and third columns.  In the forth column
we have indicated the involution $\theta$ that defines the embedding
of $K$ in $G$.

The groups $GL(n,\bbR)$, $GL(n,\bbC)$ and $GL(n,\bbH)$ denote the
groups of $n \times n$ invertible matrices with entries from the set
of real numbers, complex numbers, and quaternions respectively.  We
denote the indefinite unitary groups by $U(p,q)$ as usual.  Next we
have:
\[
\begin{array}{|c||c|c|c|c|} \hline
   G_0     &                  G                  & \twoline{\text{Defining}}{\text{Form}}  &                 K                  &                \theta                        \\ \hline \hline
 SO^*(2m)   &             SO_{2m}(\bbC)          &     D_m                                 &             GL_m(\bbC)             &  g \mapsto   -J_m      (g^{-1})^t  J_m       \\ \hline
Sp(m,\bbR) &            Sp_{2m}(\bbC)            &     J_m                                 &             GL_m(\bbC)             &  g \mapsto    D_m      (g^{-1})^t  D_m       \\ \hline
  O(p,q)   &            O_{p+q}(\bbC)            &     I_{p+q}                             &     O_p(\bbC) \times O_q(\bbC)     &  g \mapsto    I_{p,q}   g          I_{p,q}   \\ \hline
 Sp(p,q)   &          Sp_{2(p+q)}(\bbC)          &     C_{p+q}                             & Sp_{2p}(\bbC) \times Sp_{2q}(\bbC) &  g \mapsto    I_{2p,2q} g          I_{2p,2q} \\ \hline
\end{array}
\] where
\[
    D_m =
    \left[ \begin{array}{cc}   0 & I_m \\ I_m &   0 \end{array} \right],
    \text{ and }
    C_m =
\left[%
\begin{array}{cccc}
 J_1      &  0      & \cdots  & 0       \\
  0       & J_1     & \cdots  & 0       \\
  \vdots  &  \vdots & \ddots  & \vdots  \\
  0       &  0      & \cdots  & J_1      \\
\end{array}%
\right] (\text{$m$ copies of $J_1$}).
\] (Recall: $J_1 = \left[ \begin{array}{cc}   0 & 1 \\ -1 & 0 \end{array}
\right]$.) The defining form for $G$ is $(x,y) = x M y^t$ with $M$ a
square matrix denoted in the table.

Note the following notation convention in the case of the symplectic
group.  The subscript $n$ of $Sp_n(\bbC)$ refers to the dimension of
the defining representation, while the $m$ in $Sp(m)$ refers to the
rank (ie: dimension of a maximal torus).  We have $n = 2m$.  We use
this convention since it is used in \cite{GW98}.  In our notation,
$Sp(m,\bbR)$ will denote split real form of $Sp_{2m}(\bbC)$.  We
denote the compact form as $Sp(m)$ (without the $\bbR$).

Among the four groups in the above table are $Sp(m,\bbR)$, the
indefinite orthogonal, $O(p,q)$ and symplectic groups, $Sp(p,q)$,
and the group $SO^*(2m)$, which is defined as:
\[
    SO^*(2m) = \left\{ g \in SO_{2m}(\bbC)| \; -J_m\bar{g}J_m = g \right\}
\]
where $g \mapsto \bar{g}$ denotes complex conjugation of the matrix
entries of $g$. Note that this definition (as well as all the
others) can be found in a new version of Chapter 1 of \cite{GW98}
(See: \url{http://www.math.ucsd.edu/~nwallach/newch1-g-w.pdf}).

Lastly, we will consider:
\[
\begin{array}{|c||c|c|} \hline
   G_0     &                  G                  &          K         \\ \hline \hline
GL(n,\bbC) &     GL_n(\bbC) \times GL_n(\bbC)    &      GL_n(\bbC)    \\ \hline
O(n,\bbC)  &      O_n(\bbC) \times O_n(\bbC)     &      O_n(\bbC)     \\ \hline
Sp(m,\bbC) &  Sp_{2m}(\bbC) \times Sp_{2m}(\bbC) &    Sp_{2m}(\bbC)   \\ \hline
\end{array}
\]  where the involution in all three cases is $(g_1, g_2) \mapsto
(g_2,g_1)$ with $(g_1,g_2) \in G$.

The complex groups $GL_n(\bbC)$, $O_n(\bbC)$ and $Sp_{2m}(\bbC)$ are
simultaneously linear algebraic groups and real Lie groups.  The
fact that they may be viewed in these two situations could cause a
certain ambiguity which we resolve as follows:  When we write
$GL_n(\bbC)$ (resp. $O_n(\bbC)$, $Sp_{2m}(\bbC)$) we refer to the
complex algebraic group, while $GL(n,\bbC)$ (resp. $O(n,\bbC)$,
$Sp(m,\bbC)$) will refer to the real Lie group.  Certainly they are
the same as an abstract groups, but the paper is more easily
understood if we keep track of the categories that contain these
objects.

To summarize, we will refer to these ten real forms of the classical
groups.   In seven non-complex cases, the complex conjugation defining
$G_0$ may be taken as $\tau(g) = \theta\left( (\overline g^t)^{-1} \right)$,
while in the three complex cases
$\tau(g_1,g_2) =((\overline g_2^t)^{-1},(\overline g_1^t)^{-1})$.

\subsection{Hilbert Series.}\label{subsec_hilbert_series} In Section
\ref{subsec_stability} it will be convenient to state an
interpretation of the main result in terms of Hilbert series, so we
introduce the following notation:
\begin{dfn*}
        Given a symmetric pair ($G,K$), let $\fg = Lie(G)$ and set:
    \begin{eqnarray*}
        \H(\fg,K; t) &:=& \sum_{d=0}^\infty h_d(\fg, K) t^d
        \text{ where: }    h_d(\fg,K) := \dim \left[\mcS^d(\fg)\right]^K
            \text{ (for $d \in \bbN$)}.
    \end{eqnarray*}
    When convenient we will also write $\H(G_0; t)$ for $\H(\fg, K;t)$ where $G_0$ is
    the corresponding real form.
\end{dfn*}

For any one of the ten families addressed in this paper, it is a
consequence of the main theorem (\ref{subsec_main_thm}) that for
fixed $d$ the numbers $h_d(\fg,K)$ stabilize for sufficiently large
defining parameters.  We set up some notation to more plainly
indicate this phenomenon.

In the case of (say) $GL(n,\bbC)$ it is convenient to write:
\[
    \lim \H(GL(n,\bbC);t) :=
    \sum_{d=0}^\infty \left[ \lim_{n \rightarrow \infty} h_d(GL(n,\bbC)) \right] t^d
\] as the term-by-term limit.  We will see that these limits exist from
Equation \ref{eq_main_thm_GL(n,C)}.  The resulting power series
obtained in this way will be referred to as the ``stable limit". We
use the same (analogous) notation for the other groups. In the cases
where the group $G_0$ is indexed by two parameters (ie: $U(p,q)$,
$O(p,q)$, $Sp(p,q)$) the stable limit corresponds to $\min(p,q)
\rightarrow \infty$.

Using this notation, we can express the following corollary to the
main theorem that we explain in Section \ref{subsec_stability}:
\[
\lim \H \left( GL(n,\bbR) \right) = \lim \H \left( Sp(m,\bbR)
\right) = \lim \H \left( GL(n,\bbH) \right) = \lim \H \left(
SO^*(2m)\right)
\]
\[
\lim \H \left( O(n,\bbC) \right) = \lim \H \left( Sp(p,q)\right) =
\lim \H \left( Sp(m,\bbC) \right) = \lim \H \left( O(p,q)\right)
\]
\[
\lim \H \left( GL(n,\bbC) \right) = \lim \H \left( U(p,q)\right)
\]
Establishing the above equalities are one of the purposes of this
paper.

\subsection{Notation for partitions.}\label{subsec_combinatorics}

A partition, $\lambda$, with $k$ parts, is a positive integer
sequence $\lambda_1 \geq \lambda_2 \geq \ldots \geq \lambda_k > 0$.
The number of parts of a partition is called the \emph{length} of
the partition, while the sum of the parts is called the \emph{size}
of the partition.  In general, we use the same notation for
partitions as is done in standard references such as \cite{Macd95,
Stan97, Stan99}. For example, we write $\ell(\la)$ to denote the
length (or number of parts) of a partition, $|\la|$ for the size of
a partition (i.e., $|\la| = \sum_i \la_i$).  If $\la$ has size $m$
then we will write $\la \vdash m$.  Also, $\la^\prime$ denotes the
transpose (or conjugate) of $\la$ (i.e., $(\la^\prime)_i = |\{
\la_j: \la_j \geq i\}|$). Note that $\ell(\la^\prime) = \la_1$.

We implicitly identify a partition with its Young diagram.  Thus,
many terms (such as the conjugate partition) have clear meanings.
For example, we will say that $\la$ has \emph{even rows} if it is of
the form: $2\delta_1 \geq 2\delta_2 \geq \ldots \geq 2\delta_k$,
which we denote by: $\la = 2\delta$.  In the same way, we say that
$\la$ has \emph{even columns} if $\la = (2\delta)^\prime$ for some
partition $\delta$.  These partitions play a fundamental role for
us.

It will be very important to keep track of the length of the
partitions that arise in this paper.  We shall see that one
innocent, but key fact is that $\ell(\la) \leq |\la|$.  Also, if
$\la$ is of the form $2\delta$ then $\ell(\la) = \ell(\delta)$,
while $|\la| = 2|\delta|$.  Thus $\ell(\la) \leq 2|\delta|$.  On the
other hand, if $\la$ is of the form $(2\delta)^\prime$ then
$\ell(\la) \leq 2|\delta|$, but we cannot say anything about the
relationship between $\ell(\la)$ and $\ell(\delta)$.

If $\la$ denotes a partition of length $k$ and $n \geq k$, then
$\la$ defines an $n$-tuple of non-negative integers by padding out
with zeros on the right end (ie:  $(\la_1, \la_2,$ $ \cdots, $
$\la_k, 0, 0, \cdots, 0)$).  In this way, we have an injection of
the set of partitions of length at most $n$ into the lattice
$\bbZ^n$.

Another, slightly non-standard notation that we will use involves
embedding \emph{pairs} of partitions into $\bbZ^n$.  For a positive
integer $n$, let $(\mu,\nu)$ denote the $n$-tuple determined by
partitions $\mu$ and $\nu$ by:
\begin{equation}\label{eq_double_partition_notation}
    \begin{array}{ccc}
        (\mu,\nu) &:=& \underbrace{
        \left(
        \mu_1, \mu_2, \cdots, \mu_p, 0, \cdots, 0, -\nu_q, \cdots, -\nu_1
        \right)
        } \\ & & n
    \end{array}
\end{equation}
where we assume that $\ell(\mu) = p$, $\ell(\nu) = q$ and $p+q \leq
n$.  It is clear that every weakly decreasing $n$-tuple of integers
corresponds uniquely to an ordered pair of partitions in this way.

\subsection{Notation for the irreducible finite dimensional
representations.}\label{subsec_irreps}

We now fix a notation for irreducible finite dimensional
representations of $GL_n(\bbC)$.  To be precise, we will view these
groups in the category of linear algebraic groups over $\bbC$ with
morphisms being homomorphisms of the groups which are regular maps
of the underlying affine varieties.  In particular, we consider
only regular representations.\footnote{In other sources, a \emph{regular representation} is referred to as a \emph{rational} representation, but with the same meaning as we indicate here.}

For connected groups, such representations are parameterized by
dominant integral characters of a maximal torus using the theorem of
the highest weight. In the classical cases, weights are
parameterized in the standard coordinates (see \cite{GW98}) using
$n$-tuples of integers.

If a representation of $GL_n(\bbC)$ has polynomial matrix
coefficients then the $n$-tuple indexing the highest weight of the
representation has non-negative integer components. Therefore, we
may index this $n$-tuple by a partition $\la$ such that $\ell(\la)
\leq n$.

In general, the matrix coefficients of a representation of
$GL_n(\bbC)$ are not polynomials, but rather, are rational
functions in the matrix entries.  For these representations, the highest weight is indexed
by an $n$-tuple with negative components.  In light of
this, we set up the following notation which will be used in Section
\ref{sec_closed_expressions}.

Given non-negative integers $p$ and $q$ such that $n \geq p+q$ and
non-negative integer partitions $\a$ and $\b$ with $p$ and $q$ parts
respectively, let $\F{(\a,\b)}{n}$ denote the irreducible rational
representation of $GL_n$ with highest weight given by the $n$-tuple
$(\a,\b)$ (see Equation \ref{eq_double_partition_notation}). If $\b
= (0)$ (as is the case for polynomial representations) then we will
write $\F{\a}{n}$ for $\F{(\a, \b)}{n}$.  This will be our notation
for representations with polynomial matrix coefficients.

Note that if $\a = (0)$
then $\dual{\F{\b}{n}}$ is equivalent to $\F{(\a, \b)}{n}$.  In
general, the dual of the representation corresponding to $(\a,\b)$
has a highest weight which is given by the $n$-tuple $(\b, \a)$ in
the standard coordinates.

If we choose the standard basis for $\bbC^n$ we obtain a
representation of $GL_n(\bbC)$ given by the identity map $GL_n(\bbC)
\rightarrow GL_n(\bbC)$, which we will call the \emph{standard
representation}, and is denoted $\F{((1),(0))}{n}$ in our notation.
We will refer to ``the" standard representation later on in the
paper -- this is what we mean.

\bigskip

If $G = GL_n(\bbC) \times GL_m(\bbC)$ then representatives of the
irreducible representations of $G$ are of the form
$\F{(\mu_1,\mu_2)}{n} \tot \F{(\nu_1,\nu_2)}{m}$.  Note that we use
the symbol $\tot$ to denote the ``outer tensor product", which we
view as a representation of $G$.  If $k:=n=m$ then we may restrict
this representation to the diagonal $GL_k(\bbC)$; we denote this
representation $\F{(\mu_1,\mu_2)}{k} \ot \F{(\nu_1,\nu_2)}{k}$ where
the $\ot$ denotes the ``inner" tensor product of the
$GL_k(\bbC)$-representations.

\subsection{Notation for multiplicities}\label{subsec_multiplicities}
Given completely reducible representations, $V_1$ and $V_2$ of
complex algebraic groups $G_1$ and $G_2$ respectively, together with
an embedding $G_1 \hookrightarrow G_2$, we let
\[
\left[V_1,V_2\right] = \dim \text{Hom}_{G_1} \left(V_1,V_2 \right)
\]
where $V_2$ is regarded as a representation of $G_1$ by restriction.
If $V_1$ is irreducible, then $[V_1,V_2]$ is the multiplicity of
$V_1$ in $V_2$.  This cardinal may of course be infinite if $V_1$ or
$V_2$ is infinite dimensional.  For much of the paper, the
restriction to $G_1$ will be implicit, but we mention it here to be
precise.

If $G$ is a reductive algebraic group over $\bbC$, let $(V^\la)_{\la
\in \widehat G}$ denote representatives of irreducible regular
representations of $G$.  If $V$ is a completely reducible
representation of $G$, then set $m_\la := [V^\la, V]$.  In this
paper, we will always have $m_\la < \infty$, and we will indicate
the decomposition of $V$ into irreducible representations of $G$
(with multiplicity) by the expression:
\[
    V \cong \bigoplus m_\la V^\la,
\] which of course is shorthand for:
\[
    \begin{array}{rccc}
        V &\cong& \bigoplus_{\la \in \widehat G} & \underbrace{V^\la
        \oplus V^\la \oplus \cdots \oplus V^\la} \\ & & & m_\la \text{ copies. }
    \end{array}
\]

\subsection{Littlewood-Richardson coefficients.}\label{subsec_LRC}

In \cite{Macd95} the Littlewood-Richardson coefficients are defined
as the structure constants for multiplication in the ring of
symmetric polynomials at the Schur basis.  That is,
\[
    s_\mu \; s_\nu = \sum_\la \c{\la}{\mu}{\nu} \, s_\la
\] where $s_\g$ denotes the Schur function indexed by the
partition $\g$.

For most of \cite{Macd95} one works in a ring having infinitely many
variables, so there is no restriction on the number of parts of
$\la$, $\mu$ and $\nu$.  If one passes to finitely many variables
the only change that needs to be made is to make sure that all
partitions involved do not have more parts than the number of
variables.  Keeping track of this caveat, one can interpret the
Schur function as the character of a representation of $GL_n(\bbC)$.
We obtain from this interpretation:

\begin{prop}\label{prop_tensor_LRC}
If we regard the irreducible $GL_n(\bbC)$ $\times$ $GL_n(\bbC)$
representation $\F{\mu}{n}$ $\tot$ $\F{\nu}{n}$ as a representation
of the diagonally embedded $GL_n(\bbC)$ by restriction then the
decomposition into irreducibles is given by:
\[
    \F{\mu}{n} \ot \F{\nu}{n} \cong \bigoplus \c{\la}{\mu}{\nu} \F{\la}{n}
\]
where the sum is over all partitions $\la$ such that $|\la| = |\mu| + |\nu|$ and $\ell(\la) \leq n$.
\end{prop}

In \cite{GW98} and \cite{Howe95} the following standard result is
proved, which we will need in the proof of the main theorem.

\begin{prop}\label{prop_restrict_LRC}  If we regard the irreducible
$GL_{n+m}(\bbC)$-representation as a $GL_n(\bbC) \times
GL_m(\bbC)$-representation by restriction under the embedding:
\[
    \left[
    \begin{array}{cc}
        GL_n(\bbC) & 0         \\
        0          & GL_m(\bbC)
    \end{array}
    \right] \subset GL_{n+m}(\bbC),
\]
then we have:
\[
        \F{\la}{n+m} \cong \bigoplus_{\mu, \nu}
        \c{\la}{\mu}{\nu} \F{\mu}{n} \tot \F{\nu}{m},
\]
where the sum is over all partitions $\mu$ and $\nu$ such that
$\ell(\mu) \leq n$ and $\ell(\nu) \leq m$.
\end{prop}

\bigskip

\subsection{Main Theorem.}\label{subsec_main_thm} We proceed in cases:
\begin{description}
\item[CASE $G_0 = U(p,q)$ ]$\empty$ \\
  For $\fg = \gl_{p+q}(\bbC)$
  and $K = GL_p(\bbC) \times GL_q(\bbC)$,
  \begin{equation}\label{eq_main_thm_U(p,q)}
    h_d(\fg, K) = \sum (\c{\la}{\mu}{\nu})^2
  \end{equation}
  where the sum is over all partitions $\la$, $\mu$, and $\nu$ such that:
  $|\la| = |\mu| + |\nu| = d$, $\ell(\la) \leq p+q$, $\ell(\mu) \leq p$
  and $\ell(\nu) \leq q$.

\item[CASE $G_0 = GL(n, \bbR)$    ]$\empty$ \\
  For $\fg = \gl_n(\bbC)$
  and $K = O_n(\bbC)$,
  \begin{equation}\label{eq_main_thm_GL(n,R)}
    h_d(\fg, K) = \sum \c{2\la}{\mu}{\mu}
  \end{equation}
  where the sum is over all partitions $\la$, $\mu$ such that:
  $|\la| = |\mu| = d$ and $\ell(\la), \ell(\mu) \leq n$.

\item[CASE $G_0 = GL(n, \bbC)$    ]$\empty$ \\
  For $\fg = \gl_n(\bbC)$
  and $K = GL_n(\bbC)$,
  \begin{equation}\label{eq_main_thm_GL(n,C)}
    h_d(\fg, K) = \sum (\c{\la}{\mu}{\nu})^2
  \end{equation}
  where the sum is over all partitions $\la$, $\mu$, and $\nu$
  such that: $|\la| = |\mu| + |\nu| = d$,
  $\ell(\la), \ell(\mu), \ell(\nu) \leq n$.

\item[CASE $G_0 = GL(m, \bbH)$    ]$\empty$ \\
  For $\fg = \gl_{2m}(\bbC)$
  and $K = Sp_{2m}(\bbC)$,
  \begin{equation}\label{eq_main_thm_GL(n,H)}
    h_d(\fg, K) = \sum \c{(2\la)^\prime}{\mu}{\mu}
  \end{equation}
  where the sum is over all partitions $\la$, $\mu$ such that:
  $|\la| = |\mu| = d$ and $\ell((2\la)^\prime), \ell(\mu) \leq 2m$.

\item[CASE $G_0 = SO^*(2m)$      ]$\empty$ \\
  For $\fg = \so_{2m}(\bbC)$
  and $K = GL_m(\bbC)$,
  \begin{equation}\label{eq_main_thm_SO^*(2m)}
    h_d(\fg, K) = \sum \c{(2\la)^\prime}{\mu}{\mu}
  \end{equation}
  where the sum is over all partitions $\la$, $\mu$
  such that: $|\la| = |\mu| = d$ and $\ell((2\la)^\prime) \leq 2m$,
  and $\ell(\mu) \leq m$.

\item[CASE $G_0 = Sp(m, \bbR)$ ]$\empty$ \\
  For $\fg = \sp_n(\bbC)$
  and $K = GL_m(\bbC)$,
  \begin{equation}\label{eq_main_thm_Sp(m,R)}
    h_d(\fg, K) = \sum \c{2\la}{\mu}{\mu}
  \end{equation}
  where the sum is over all partitions $\la$, $\mu$
  such that: $|\la| = |\mu| = d$ and $\ell(\la) \leq 2m$, and $\ell(\mu) \leq m$.

\item[CASE $G_0 = O(n, \bbC)$     ]$\empty$ \\
  For $\fg = \so_n(\bbC) \oplus \so_n(\bbC)$
  and $K = O_n(\bbC)$,
  \begin{equation}\label{eq_main_thm_O(n,C)}
    h_d(\fg, K) = \sum \c{2\la}{(2\mu)^\prime}{(2\mu)^\prime}
  \end{equation}
  where the sum is over all partitions $\la$, $\mu$, and $\nu$
  such that: $|\la| = |\mu|+|\nu| = d$ and
  $\ell(\la)$, $\ell((2\mu)^\prime)$, $\ell((2\nu)^\prime)$ $\leq n$.

\item[CASE $G_0 = Sp(m,\bbC)$ ]$\empty$ \\
  For $\fg = \sp_{2m}(\bbC) \oplus \sp_{2m}(\bbC)$
  and $K = Sp_{2m}(\bbC)$,
  \begin{equation}\label{eq_main_thm_Sp(m,C)}
    h_d(\fg, K) = \sum \c{(2\la)^\prime}{2\mu}{2\mu}
  \end{equation}
  where the sum is over all partitions $\la$, $\mu$
  such that: $|\la| = |\mu|+|\nu| = d$ and
  $\ell((2\la)^\prime), \ell(\mu), \ell(\nu) \leq 2m$.

\item[CASE $G_0 = O(p,q)$       ]$\empty$ \\
  For $\fg = \so_{p+q}(\bbC)$
  and $K = O_p(\bbC) \times O_q(\bbC)$,
  \begin{equation}\label{eq_main_thm_O(p,q)}
    h_d(\fg, K) = \sum \c{(2\la)^\prime}{2\mu}{2\mu}
  \end{equation}
  where the sum is over all partitions $\la$, $\mu$, and $\nu$
  such that: $|\la| = |\mu|+|\nu| = d$ and
  $\ell((2\la)^\prime) \leq p+q$, $\ell(\mu) \leq p$,  $\ell(\nu) \leq q$.

\item[CASE $G_0 = Sp(p,q)$      ] $\empty$ \\
  For $\fg = \sp_{2p+2q}(\bbC)$
  and $K = Sp_{2p}(\bbC) \times Sp_{2q}(\bbC)$,
  \begin{equation}\label{eq_main_thm_Sp(p,q)}
    h_d(\fg, K) = \sum \c{2\la}{(2\mu)^\prime}{(2\mu)^\prime}
  \end{equation}
  where the sum is over all partitions $\la$, $\mu$, and $\nu$
  such that: $|\la| = |\mu|+|\nu| = d$ and
  $\ell(\la) \leq 2(p+q)$, $\ell((2\mu)^\prime) \leq 2p$,
  $\ell((2\nu)^\prime) \leq 2q$.
\end{description}
\begin{proof}
See Section \ref{sec_main_thm_proof} for a case--by--case analysis.
\end{proof}

\begin{rem}
Note that the cases $GL(n, \bbR)$ and $GL(n, \bbH)$ are also
addressed in \cite{Will00} and then more thoroughly in
\cite{Will01}.
\end{rem}

\subsection{Stability properties.}\label{subsec_stability}  In each
of the ten cases addressed in the main theorem the summation range
involves a condition on both the \emph{size} and the \emph{length}
of the partitions involved.  The parameter describing the bound for
the length is independent of the parameter for the size.  For fixed
size, as the length parameter goes to infinity the condition on the
length is automatic due to the fact that for a partition $\la$,
$\ell(\la) \leq |\la|$. This innocent combinatorial fact allows us
to conclude the following stable behavior for a fixed $d$:
\begin{itemize}
\item Examining Equations \ref{eq_main_thm_GL(n,R)} and
\ref{eq_main_thm_Sp(m,R)} we see that:
\begin{equation}\label{eq_stable_GL(n,R)-Sp(m,R)}
    h_d\left( GL(n,\bbR) \right) = h_d\left( Sp(m,\bbR) \right)
\end{equation}
provided $d \leq \min(n,m)$;
\item Examining Equations \ref{eq_main_thm_GL(n,H)} and
\ref{eq_main_thm_SO^*(2m)} we see that:
\begin{equation}\label{eq_stable_GL(n,H)-O*(2m)}
    h_d\left( GL(n,\bbH) \right) = h_d\left( SO^*(2m) \right)
\end{equation}
provided $d \leq \min(n,m)$;
\item Examining Equations \ref{eq_main_thm_O(n,C)} and
\ref{eq_main_thm_Sp(p,q)} we see that:
\begin{equation}\label{eq_stable_O(n,C)-Sp(p,q)}
    h_d\left( O(n,\bbC) \right) = h_d\left( Sp(p,q) \right)
\end{equation}
provided $d \leq \min(\frac{n}{2}, p,q)$;
\item Examining Equations \ref{eq_main_thm_Sp(m,C)} and
\ref{eq_main_thm_O(p,q)} we see that:
\begin{equation}\label{eq_stable_Sp(m,C)-O(p,q)}
    h_d\left( Sp(m,\bbC) \right) = h_d\left( O(p,q) \right)
\end{equation}
provided $d \leq \min(m, \frac{p}{2},\frac{q}{2})$;
\item Examining Equations \ref{eq_main_thm_GL(n,C)} and
\ref{eq_main_thm_U(p,q)} we see that:
\begin{equation}\label{eq_stable_GL(n,C)-U(p,q)}
    h_d\left( GL(n,\bbC) \right) = h_d\left( U(p,q) \right)
\end{equation}
provided $d \leq \min(n,p,q)$.
\end{itemize}

A direct consequence of the above equalities are the following
stable limits:
\begin{eqnarray}
    \lim \H \left( GL(n,\bbR) \right) &=& \lim \H \left( Sp(m,\bbR)
    \right) \\
    \lim \H \left( GL(n,\bbC) \right) &=& \lim \H \left( U(p,q)
    \right) \\
    \lim \H \left( GL(n,\bbH) \right) &=& \lim \H \left( SO^*(2m)
    \right) \\
    \lim \H \left( O(n,\bbC) \right) &=& \lim \H \left( Sp(p,q)
    \right) \\
    \lim \H \left( Sp(m,\bbC) \right) &=& \lim \H \left( O(p,q)
    \right)
\end{eqnarray}
Certainly, one might expect that these formulas follow for other
arguments involving, for example, the theory of dual pairs (see
\cite{Howe89-1, Howe89-2}). Such an approach would indeed be very
interesting.

From the point of view presented here, we can exploit combinatorial
facts concerning the Littlewood-Richardson coefficients.  For
example, in \cite{Macd95} and \cite{Stan99} the following well know
fact is presented:
\begin{prop}  For any partitions $\alpha$. $\beta$, $\gamma$ we have:
$\c{\gamma}{\alpha}{\beta} =
\c{\gamma^\prime}{\alpha^\prime}{\beta^\prime}$.
\end{prop}

Consequently, we obtain that the sums in Equations
\ref{eq_main_thm_GL(n,R)}, \ref{eq_main_thm_GL(n,H)},
\ref{eq_main_thm_SO^*(2m)}, and \ref{eq_main_thm_Sp(m,R)} are all
equal for sufficiently large defining parameters. Moreover, we see
another equality between Equations \ref{eq_main_thm_O(n,C)},
\ref{eq_main_thm_Sp(m,C)}, \ref{eq_main_thm_O(p,q)},
\ref{eq_main_thm_Sp(p,q)}.  Thus, the ten cases represent only
three:
\[
\lim \H \left( GL(n,\bbR) \right) = \lim \H \left( Sp(m,\bbR)
\right) = \lim \H \left( GL(n,\bbH) \right) = \lim \H \left(
SO^*(2m)\right)
\]
\[
\lim \H \left( O(n,\bbC) \right) = \lim \H \left( Sp(p,q)\right) =
\lim \H \left( Sp(m,\bbC) \right) = \lim \H \left( O(p,q)\right)
\]
\[
\lim \H \left( GL(n,\bbC) \right) = \lim \H \left( U(p,q)\right)
\]

\subsection{Acknowledgements}

The author would like to thank several people for conversations
related to the results in this paper including:  Thomas Enright
(UCSD), Roger Howe (Yale), Markus Hunziker (Baylor), Soo Teck Lee
(NUS), Richard Stanley (MIT), Eng Chye Tan (NUS), and Nolan Wallach
(UCSD).

Also, the material in this paper was the subject of a lecture given
in the \emph{Geometry, Symmetry and Physics} seminar at Yale
University organized by Igor Frenkel and Gregg Zuckerman in the
Spring of 2003. It was also the topic of a talk at the Special
Session on Recent Advances in Non-Commutative Algebra in the
Southeastern Sectional Meeting of the AMS in Western Kentucky
University during March 2005, organized by Ellen Kirkman of Wake
Forest University.

\bigskip
This research was supported by NSA Grant \#H98230-05-1-0078.

\section{Background for the proofs.}\label{sec_background}
Interestingly, we require relatively little machinery to prove the
main theorem.  In fact, the only results required are multiplicity
free spaces for $GL_n(\bbC)$ and (related) standard results concerning the
Littlewood-Richardson coefficients.  We state the necessary results
here for reference.
\subsection{The Cartan-Helgason theorem.}
One of the tools that we will use for the proof of the main theorem
is the Cartan-Helgason theorem (see \cite{GW98} or \cite{Howe95}).
This fact amounts to the assertion that if $(G,K)$ is a symmetric
pair then the $G$ decomposition of regular functions on $G/K$ is
multiplicity free. The irreducible regular $G$-representations that
occur in $\mcO(G/K)$ are precisely those representations with a
$K$-invariant vector.  We describe here four classical instances of
this theorem, the first of which we describe as a remark since it is
a consequence of Schur's lemma:
\begin{rem}\label{rem_Cartan_Helgason_GxG}
For any reductive group $K$, we may symmetrically embedded $K$ in $G
:= K \times K$ as $\left\{(k,k)|k \in K \right\}$. The irreducible
regular representations of $G$ may be taken as $V_1 \tot V_2$ where
$(\rho_1,V_1)$ and $(\rho_2,V_2)$ are irreducible representations of
$K$.  We then obtain:
\[
    \dim \left( V_1 \otimes V_2 \right)^K =
    \left\{%
\begin{array}{ll}
    1, & \hbox{if $V_1 \cong V_2^*$;} \\
    0, & \hbox{otherwise,} \\
\end{array}%
\right.
\]
by noting that $\text{Hom}_K (V_1, V_2) \cong$ $\left(V_1^* \otimes
V_2 \right)^K$ and invoking Schur's lemma.
\end{rem}
\subsubsection{Symmetric pair: $(GL_n(\bbC), O_n(\bbC))$}\label{thm_Cartan-Helgason_GLnOn}
An irreducible representation of $GL_n(\bbC)$ may be regarded as a representation of $O_n(\bbC)$
by restriction.  In so doing, we have the following
\begin{thm*}
For an integer partition $\la$ with $\ell(\la) \leq n$:
\[
    \dim \left( \F{\la}{n} \right)^{O_n(\bbC)} =
    \left\{%
\begin{array}{ll}
    1, & \hbox{$\la_i$ even for all $i$;} \\
    0, & \hbox{otherwise.} \\
\end{array}%
\right.
\]
\end{thm*}
\subsubsection{Symmetric pair: $(GL_{2m}(\bbC), Sp_{2m}(\bbC))$}\label{thm_Cartan-Helgason_GL2mSp2m}
As in the last case, an irreducible representation of $GL_{2m}(\bbC)$ may be regarded as a representation of $Sp_{2m}(\bbC)$
by restriction.  We have the following
\begin{thm*}
For an integer partition $\la$ with $\ell(\la) \leq 2m$:
\[
    \dim \left( \F{\la}{2m} \right)^{Sp_{2m}(\bbC)} =
    \left\{%
\begin{array}{ll}
    1, & \hbox{$(\la^\prime)_i$ even for all $i$;} \\
    0, & \hbox{otherwise.} \\
\end{array}%
\right.
\]
\end{thm*}
\subsubsection{Symmetric pair: $(GL_{p+q}(\bbC), GL_p(\bbC) \times GL_q(\bbC))$}\label{thm_Cartan-Helgason_GLnGLpGLq}
As in the last cases, an irreducible representation of $GL_{p+q}(\bbC)$ may be regarded as a representation of $GL_p(\bbC) \times GL_q(\bbC)$
by restriction where we embed $GL_p(\bbC) \times GL_q(\bbC)$ in $GL_{p+q}(\bbC)$ as:
\[
    \left[
    \begin{array}{cc}
        GL_p(\bbC) &          0 \\
        0          & GL_q(\bbC)
    \end{array}
    \right] \subset GL_{p+q}(\bbC).
\]
We have the following
\begin{thm*}
  Given $n = p+q$, let $\mu$ and $\nu$ be partitions with at most $p$ and $q$ parts
respectively.  Then,
\[
\dim( \F{(\mu, \nu)}{n} )^{GL_p \times GL_q} =
\begin{cases}
    1 & \text{if $\mu = \nu$}, \\
    0 & \text{otherwise}.
\end{cases}
\]
\end{thm*}

\subsection{Classical Multiplicity Free Spaces for $GL_n(\bbC)$.}
We recall, in our notation, three celebrated multiplicity free
representations for $GL_n(\bbC)$.  See for \cite{GW98} or
\cite{Howe95} for proofs.

\subsubsection{The $k \times m$ Matrices.}
$\empty$\\
We define the following action of $GL_k(\bbC)\times
GL_m(\bbC)$ on the algebra of complex valued polynomial function on
the $k \times m$ complex matrices. For $(g,h) \in GL_k(\bbC) \times
GL_m(\bbC)$, $X \in M_{k,m}$ and $f \in \mcP(M_{k,m})$, set: $(g,h)
\cdot f (X) = f(g^t X h)$. Under this action we can identify:
$\mcP(M_{k,m}) \cong \mcS(\bbC^k \tot \bbC^m)$ as a representation
of the group $GL_k(\bbC) \times GL_m(\bbC)$. For our purposes, the
decomposition into irreducible representations is of particular
interest:
\begin{thm}\label{thm_Mpq} As a $GL_k(\bbC) \times GL_m(\bbC)$-representation,
\[  \mcS^d(\bbC^k \tot \bbC^m) \cong
    \bigoplus_{\twoline{\la:|\la|=d}{\ell(\la)\leq \min(k,m)}} \F{\la}{n} \tot \F{\la}{m}.
\]
\end{thm}

We also state the following dual version of the above theorem.
Define the action of $GL_k(\bbC)\times GL_m(\bbC)$ on $\mcP(M_{k,m})$ by: for $(g,h) \in GL_k(\bbC) \times GL_m(\bbC)$, $X \in M_{k,m}$ and $f \in \mcP(M_{k,m})$, set: $(g,h) \cdot f (X) = f(g^{-1} X h)$.
Under this action we can identify: \\
$\mcP(M_{k,m}) \cong \mcS(\dual{\bbC^k} \tot \bbC^m)$ as a representation of the group $GL_k(\bbC) \times GL_m(\bbC)$.

\begin{cor}\label{cor_dual_Mpq}As a $GL_k(\bbC) \times GL_m(\bbC)$-representation,
\[  \mcS \left(\dual{\bbC^k} \tot \bbC^m \right) \cong
    \bigoplus_{\twoline{\la:|\la|=d}{\ell(\la)\leq \min(k,m)}}
    \dual{\F{\la}{k}} \tot \F{\la}{m}.
\]
\end{cor}

\subsubsection{The Symmetric Matrices.}\label{thm_SMn}
$\empty$\\
We define the following action of $GL_n(\bbC)$ on the
algebra of complex valued polynomial function on the symmetric $n
\times n$ complex matrices, $SM_n$. For $g \in GL_n(\bbC)$, $X \in
SM_n$ and $f \in \mcP(SM_n)$, set: $g \cdot f (X) = f(g^t X g)$.
Under this action we can identify: $\mcP(SM_n) \cong \mcS(\mcS^2
\bbC^n)$ as a representation of the group $GL_n(\bbC)$.  We have the
multiplicity free decomposition:
\begin{thm*} As a $GL_n(\bbC)$-representation,
\[  \mcS^d(\mcS^2 \bbC^n) \cong
    \bigoplus_{\twoline{\la:|\la|=d}{\ell(\la)\leq n}} \F{2\la}{n}.
\]
\end{thm*}

\subsubsection{The Skew-Symmetric Matrices.}\label{thm_AMn}
$\empty$\\
We define an analogous action of $GL_n(\bbC)$ on the
algebra of complex valued polynomial function on the skew-symmetric
$n \times n$ complex matrices, $AM_n$. For $g \in GL_n(\bbC)$, $X
\in AM_n$ and $f \in \mcP(AM_n)$, set: $g \cdot f (X) = f(g^t X g)$.
Under this action we can identify: $\mcP(AM_n) \cong \mcS(\wedge^2
\bbC^n)$ as a representation of the group $GL_n(\bbC)$.  We have the
multiplicity free decomposition:
\begin{thm*} As a $GL_n(\bbC)$-representation,
\[  \mcS^d(\wedge^2 \bbC^n) \cong
    \bigoplus_{\twoline{\la:|\la|=d}{\ell((2\la)^\prime)\leq n}} \F{(2\la)^\prime}{n}.
\]
\end{thm*}

\section{Proof of the Main Theorem.}\label{sec_main_thm_proof}

We now proceed to prove the main theorem.
We do this in case-by-case analysis, however each is presented in the following manner:
We start with the adjoint representation, $\fg$, of $G$, which gives rise to a representation of
$G$ on the degree $d$ symmetric tensors, $\mcS^d(\fg)$.  Our goal is to compute the dimension of
invariants when this representation is restricted to $K$,
\begin{eqnarray*}
           G    & \rightarrow GL\left( \mcS^d(\fg) \right) \\
           \cup &                       \\
           K    &
\end{eqnarray*}

This representation may be directly analyzed using branching rules from $G$
to $K$, but in general this a not the best route to the result.
Instead, it will be convenient to realize that $K$ is symmetrically embedded
in another group $\cK$ which also acts on $\mcS^d(\fg)$.  Furthermore, there exists
yet another group, $\cG$ in which both $G$ and $\cK$ are symmetrically embedded.
The picture is:
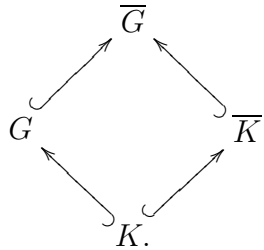
\[
\xymatrix{
                   & \cG                 &             \\
G \ar@{^{(}->}[ru] &                     & \cK \ar@{^{(}->}[lu] \\
                   & K. \ar@{^{(}->}[lu] \ar@{^{(}->}[ru]   &
}
\]
The idea is to exploit the existence of a certain graded,
$\cG$-multiplicity free representation, $\cmcS = \bigoplus \cmcS[d]$
which as a graded $G$-representation is equivalent to $\mcS(\fg)$.
This representation is of the form $\cmcS := \mcS(\cW)$ for some
$\cG$-representation $\cW$.

Understanding the graded $K$-module structure of $\mcS(\fg)$ will be
equivalent to understanding the graded $K$-module structure of
$\cmcS$.  We will see that it is technically much easier to first
decompose relative to the $\cK$ action and then use a classical
instance of the Cartan-Helgason theorem to determine the dimension
of the space of $K$-invariants.

\subsection{CASE: $G_0 = U(p,q)$}  Let $n := p+q$.
\[
\begin{array}{|rcl|rcl|} \hline
G        & = & GL_n(\bbC)   \hspace{15mm}  & K    & \cong & GL_p(\bbC) \times GL_q(\bbC) \\
\hline
\end{array}
\]
The group $K$ is symmetrically embedded in $G$ in the standard way as:
\[K := \left[
\begin{array}{cc}
 GL_p(\bbC) &     0      \\
     0      & GL_q(\bbC) \\
\end{array}
\right] \subset G
\]

We then define $\cG$ and $\cK$ as:
\[
\begin{array}{ccc}
\cG:=G \times G & \supset & \cK:=K \times K \\
\cup            &         & \cup           \\
G               & \supset & K
\end{array}
\] with the obvious vertical diagonal embeddings.

Let $W = \bbC^n$ denote the standard representation of $G$. Under
the adjoint representation of $G$ we have: $\fg  \cong W^* \ot W$.
This representation is the restriction of the irreducible
$\cG$-representation $\cW \cong W^* \tot W$. Let $\cmcS :=
\mcS(\cW)$.  By Corollary \ref{cor_dual_Mpq}, we have:
\begin{equation}
    \mcS^d \left(\cW \right) \cong   \bigoplus_{\twoline{\la: |\la|=d}{\ell(\la)\leq n}} \dual{ \F{\la}{n} } \tot \F{\la}{n}
\end{equation} under the action of $\overline G$.
We will exploit this fact by restricting to the subgroup, $\cK
\subseteq \cG$. We then apply Proposition \ref{prop_restrict_LRC} to
each occurrence of $\F{\la}{p+q}$ to obtain the following
decomposition of the $\mcS^d(\cW)$ as a $\cK$-representation,
\[
    \mcS^d(W^* \tot W) \cong \bigoplus \dual{\c{\la}{\mu}{\nu} \F{\mu}{p} \tot \F{\nu}{q}} \ot \left( \c{\la}{\alpha}{\beta} \F{\alpha}{p} \tot \F{\beta}{q}  \right)
\]
where the above sum is over $\alpha, \beta, \mu, \nu$ and $\la$ such that:
\begin{equation}\label{eq_Upq_summation_parameters}
\begin{array}{c}
\ell(\mu), \ell(\alpha) \leq p \\
\ell(\nu), \ell(\beta)  \leq q \\
\ell(\la) \leq n, \mbox{ and } |\la| = d.
\end{array}
\end{equation}

Schur's lemma (see Section \ref{rem_Cartan_Helgason_GxG}) implies
that we will have a $K$-invariant exactly when $\mu = \alpha$ and
$\nu = \beta$.  And so,

\[ \dim \mcS^d(\fg)^K = \sum (\c{\la}{\mu}{\nu})^2 \dim
\left[
    \dual{\F{\mu}{p} \tot \F{\nu}{q}} \tot \F{\mu}{p} \tot \F{\nu}{q}
\right]^K
\]
where the sum is over all partitions $\la$, $\mu$, and $\nu$ as in
Line \ref{eq_Upq_summation_parameters}.  Again by Schur's lemma we
see that: \\
$\dim \left[ \dual{\F{\mu}{p} \tot \F{\nu}{q}} \tot \F{\mu}{p} \tot
\F{\nu}{q} \right]^K = 1$. Equation \ref{eq_main_thm_U(p,q)}
follows.

\subsection{CASE: $G_0 = GL(n, \bbR)$   }
\[
\begin{array}{|rcl|rcl|} \hline
G   &=&  GL_n(\bbC) \hspace{15mm} & K   &\cong& O_n(\bbC) \\
\hline
\end{array}
\]
The group $K$ is symmetrically embedded in $G$ in the standard way as:
\[  K := \{g \in G| g^{-1} = g^t \} \subset G. \]
We then define $\cG$ and $\cK$ as:
\[
\begin{array}{cclccl}
\cG      &:=&G \times G                               & \supset & \cK      & := \{ (g,g) | g \in G \}                 \\
\cup     &  &                                         &         & \cup     &                                          \\
G^\prime &:=& \{ ((g^{-1})^t,g) | g \in G \} & \supset & K^\prime & := \{ (k,k) | k \in K \}.
\end{array}
\]
Note that the primes are indicated to describe the symmetric embeddings of $G$ and $K$ precisely.
In particular, note that we have embedded $GL_n(\bbC)$ in $\cG$ in two distinct ways:
\[GL_n(\bbC) \cong \cK := \{ (g,g)| g \in GL_n(\bbC) \} \subset \cG \] and
\[GL_n(\bbC) \cong  G  := \{ ((g^{-1})^t,g)| g \in GL_n(\bbC) \} \subset \cG. \]
As before, let $W = \bbC^n$ denote the standard representation of $GL_n(\bbC)$.

Under the adjoint representation of $G$ we have $\fg  \cong W^* \ot
W$, which is the restriction of the irreducible representation $W^*
\tot W$ of $\overline G$.  Note that under restriction to $K$, $W
\cong W^*$ since the regular representations of $O_n(\bbC)$ are
self-dual. Let $\cmcS := \mcS(\cW)$ where $\cW := W \tot W$. As a
representation of $K$, $\mcS(\fg) \cong \cmcS$. By Theorem
\ref{thm_Mpq} we have:
\[
    \mcS^d(\cW) \cong \bigoplus_{\twoline{\mu: |\mu|=d}{\ell(\mu)\leq n}} \F{\mu}{n} \tot \F{\mu}{n}
    \text{ as a $\cG$-representation.}
\]
We now restrict to the diagonal subgroup of $\cK$.  In doing this,
the irreducible representations $\F{\mu}{n} \tot \F{\mu}{n}$
decompose into irreducible representations of $\cK$ according to the
Littlewood-Richardson rule (see Proposition \ref{prop_tensor_LRC}).
For a partition $\nu$, the $GL_m(\bbC)$-irrep. $\F{\nu}{n}$ occurs
in $\mcS(W \ot W)$ with multiplicity $\c{\nu}{\mu}{\mu}$.  Thus,
\[
    \mcS^d(\fg)^K \cong \mcS^d(W \ot W)^{O_n(\bbC)} \cong \bigoplus \c{\nu}{\mu}{\mu} \left( \F{\nu}{n} \right)^{O_n(\bbC)}
\] where $\mu$ and $\nu$ are partitions with $\ell(\mu), \ell(\nu)$ $\leq n$ and $|\nu| = 2|\mu| = 2d$.
The Equation \ref{eq_main_thm_GL(n,R)} follows as a application of
an instance of the Cartan-Helgason Theorem (see Section
\ref{thm_Cartan-Helgason_GLnOn}).  That is, we obtain a
$K$-invariant for each $\nu$ of the form $2\la$ where $\la$ is a
partition with length at most $n$.

\subsection{CASE: $G_0 = GL(n, \bbC)$}
\[
\begin{array}{|rcl|rcl|} \hline
G   &=&  GL_n(\bbC) \times GL_n(\bbC)  \hspace{15mm} & K   &\cong& GL_n(\bbC)\\
\hline
\end{array}
\]
The group $K$ is diagonally embedded in $G$.  We define $\cG$ and the embedding of $G$ as follows:
\[
\begin{array}{ccl}
\cG      &:=& G \times G \cong \times_{i=1}^4 GL_n(\bbC)         \\
\cup     &  &                                   \\
G^\prime &:=& \{ (g,g,h,h)| g,h \in GL_n(\bbC) \} \cong G.
\end{array}
\]
We embed a group $\cK$ in $\cG$ as:
\[
\begin{array}{rccl}
\cG     & \supset & \cK      & := \{ (g,h,g,h) | g,h \in GL_n(\bbC) \}       \\
        &         & \cup     &                                               \\
K       & \cong   & K^\prime & := \{ (k,k,k,k) | k   \in GL_n(\bbC) \}.
\end{array}
\]
Let $W$ denote the standard representation of $GL_n(\bbC)$. The
adjoint representation of $GL_n(\bbC)$ is $\gl_n(\bbC) = W^* \ot W$,
and so the adjoint representation of $G = GL_n(\bbC)\times
GL_n(\bbC)$ is
\[
\begin{array}{cccl}
\fg & =     & (\gl_n( \bbC) \tot 1_\bbC) \oplus (1_\bbC \tot \gl_n(\bbC)) & \\
    & \cong & ((W^* \ot  W) \tot 1_\bbC) \oplus (1_\bbC \tot (W^* \ot W)) & \mbox{ as a $G$-representation;} \\
    & \cong &  (W^* \ot W) \oplus (W^* \ot W)                           & \mbox{ as a $  K$-representation}.
\end{array}
\]
where $1_\bbC$ denotes the trivial representation of $GL_n(\bbC)$.

We obtain a representation, $\cmcS := \mcS(\cW)$ of $\cG$ where:
\[
\cW := \left( W^* \tot W \tot 1_\bbC \tot 1_\bbC \right) \oplus
\left(1_\bbC \tot 1_\bbC \tot W^* \tot W \right).
\]  As in the other cases, $\cmcS \cong \mcS(\fg)$ as a graded $G$-representation.

Using Corollary \ref{cor_dual_Mpq} on each summand of $\cW$ we
obtain the following multiplicity free decomposition:
\[
\mcS^d(\cW) \cong \bigoplus \dual{\F{\mu}{n}} \tot \F{\mu}{n} \tot
\dual{\F{\nu}{n}} \tot \F{\nu}{n}
\]
where the sum is over all partitions $\mu$ and $\nu$ with $\ell(\mu), \ell(\nu) \leq n$ and $|\mu| + |\nu| = d$.
When we restrict to $\cK$ structure is obtained from the decomposing the inner tensor products,
$\F{\mu}{n} \ot \F{\nu}{n}$ and $\dual{\F{\mu}{n}} \ot \dual{\F{\nu}{n}}$.
The $\cK$-decomposition is:
\[
\begin{array}{ccc}
    \mcS^d(\cW) &\cong \bigoplus_{\mu,\nu} & \dual{ \bigoplus_\a \c{\a}{\mu}{\nu} \F{\a}{n} } \tot \left( \bigoplus_\b \c{\b}{\mu}{\nu} \F{\b}{n} \right) \\
                         &\cong \bigoplus_{\a,\b} & \left(\sum_{\mu,\nu} \c{\a}{\mu}{\nu}  \c{\b}{\mu}{\nu} \right) \left[\dual{\F{\a}{n}} \tot \F{\b}{n} \right]
\end{array}
\] where the sums are over $\a, \b$ such that $\ell(\a), \ell(\b) \leq n$ and $|\a| = |\b| = |\mu| + |\nu|$.  Now, we restrict to $K$,

\[
\dim \left[\mcS^d(\cW)\right]^K = \sum_{\a,\b} \left(\sum_{\mu,\nu}
\c{\a}{\mu}{\nu}  \c{\b}{\mu}{\nu} \right) \dim
\left[\dual{\F{\a}{n}} \tot \F{\b}{n} \right]^K
\]
We obtain $K$-invariants exactly when $\alpha = \beta$.  We denote
the common value of $\a$ and $\b$ by $\la$. Equation
\ref{eq_main_thm_GL(n,C)} follows.

\subsection{CASE: $G_0 = GL(m, \bbH)$   }  Let $n = 2m$.
\[
\begin{array}{|rcl|rcl|} \hline
G   &=&  GL_{2m}(\bbC) \hspace{15mm} & K   &=& Sp_{2m}(\bbC) \\
\hline
\end{array}
\]
The group $K$ is symmetrically embedded in $G$ in the standard way
as:
\[  K := \{g \in G| g^{-1} = - J_m g^t J_m \} \subset G \] where: $J_m$ is
as before.  We then define $\cG$ and $\cK$ as:
\[
\begin{array}{cclccl}
\cG      &:=&G \times G                         & \supset & \cK      & := \{ (g,g) | g \in G \}                 \\
\cup     &  &                                   &         & \cup     &                                          \\
G^\prime &:=& \{ ((Jg^{-1})^tJ,g) | g \in G  \} & \supset & K^\prime
& := \{ (k,k) | k \in K \}.
\end{array}
\]  As before, the primes are indicated to describe precisely the embeddings.
In particular, note that we have embedded $GL_{2m}(\bbC)$ in $\cG$
in two distinct ways.


This case is almost identical to the last, so we truncate the full
discussion. There are a few minor differences such as the use
Theorem \ref{thm_Cartan-Helgason_GL2mSp2m} instead of Theorem
\ref{thm_Cartan-Helgason_GLnOn}. As before, $\cmcS := \mcS(\cW)$,
where $\cW = W \tot W$ with $W$ the standard representation of $G$.
We obtain:
\[
    \mcS^d(\fg)^K \cong \mcS^d(W \ot W)^{Sp_{2m}(\bbC)} \cong \bigoplus \c{\nu}{\mu}{\mu} \left( \F{\nu}{2m} \right)^{Sp_{2m}(\bbC)}
\] where $\mu$ and $\nu$ are partitions with $\ell(\mu), \ell(\nu)$ $\leq 2m$ and $|\nu| = 2|\mu| = 2d$.
In the present case, we obtain a $K$-invariant exactly when $\nu =
(2\la)^\prime$, for some $\la$. Equation \ref{eq_main_thm_GL(n,H)}
follows.

\subsection{CASE: $G_0 = SO^*(2m)$     }
\[
\begin{array}{|rcl|rcl|} \hline
G   &=&  SO_{2m}(\bbC) \hspace{15mm} & K   &\cong& GL_m(\bbC) \\
\hline
\end{array}
\]
The group $K$ is embedded in $G$, as:
\[
    K \cong K^\prime := \left\{ \left[ \begin{array}{cc} g & 0 \\ 0 & (g^{-1})^t \end{array}\right]: g \in K \right\}  \subset G
\]  Note that we choose the form of $G$ using the matrix $D_m$ as in Section \ref{subsec_realforms}.

We define $\cG$ and the embedding of $G$ as follows:
\[
\begin{array}{ccl}
\cG      &:=& GL_{2m}(\bbC) \\
\cup     &  &                                    \\
G^\prime &:=& SO_{2m}(\bbC) \cong G
\end{array}
\] with the standard embedding of $SO_{2m}$.

Let $\cK_L$ and $\cK_R$ denote isomorphic copies of $GL_m(\bbC)$.
Let $W_L$ and $W_R$ denote the respective standard representations
of these groups.  We set $\cK := \cK_L \times \cK_R$, and embed in
$\cG$ as:
\[
\begin{array}{rcll}
\cG     & \supset &  \cK^\prime := \left\{ \left[ \begin{array}{cc} g & 0 \\ 0 & h \end{array}\right]: g \in \cK_L, h \in \cK_R \right\} \cong \cK \\
        &         &  \cup                                                                                                  \\
K       & \cong   &  K^\prime.
\end{array}
\]
Let $\overline W$ be the standard representation of $\cG$, which we
also regard as the standard representation of $G$.  When restricted
to $\cK$, $\overline W \cong W_L \tot 1_\bbC \oplus 1_\bbC \tot
W_R$, and when restricted to $K$, $\overline W \cong W \oplus
W^*$, where $W$ is the standard representation of $K$.

As a representation of $G$, $\fg \cong \wedge^2(\overline W)$. Let
$\cW := \wedge^2(\overline W)$ and $\cmcS := \mcS(\cW)$.  We view
the latter as a graded, representation of $\cG$. Using Theorem
\ref{thm_AMn}, we have a multiplicity free $\cG$-decomposition:
\[
    \mcS^d(\cW) \cong \bigoplus \F{(2\la)^\prime}{2m}
\]
where the sum is over all partitions $\la$ such that $|\la| = d$ and
$\ell((2\la)^\prime) \leq 2m$.

We now decompose relative to the action of $\cK$ using Proposition
\ref{prop_restrict_LRC}:
\[
    \mcS^d(\cW) \cong \bigoplus \c{(2\la)^\prime}{\mu}{\nu} \F{\mu}{m}
    \tot \F{\nu}{m}
\]
where $\mu$ and $\nu$ are partitions with $|\mu|+|\nu| = 2|\la|$ and
$\ell(\mu), \ell(\nu) \leq m$.

When we restrict to $K$ the decomposition becomes:
\[
    \mcS^d(\cW) \cong \bigoplus \c{(2\la)^\prime}{\mu}{\nu} \F{\mu}{m}
    \ot \dual{\F{\nu}{m}}.
\]
Clearly, we obtain a $K$ invariant exactly when $\mu = \nu$.
Equation \ref{eq_main_thm_SO^*(2m)} follows.

\subsection{CASE: $G_0 = Sp(m, \bbR)$}
\[
\begin{array}{|rcl|rcl|} \hline
G   &=&  Sp_{2m}(\bbC) \hspace{15mm} & K   &=& GL_m(\bbC) \\
\hline
\end{array}
\]
The group $K$ is diagonally embedded in $G$, as:
\[
    K \cong K^\prime := \left\{ \left[ \begin{array}{cc} g & 0 \\ 0 & -J_m (g^{-1})^t J_m \end{array}\right]: g \in K \right\}  \subset G
\] Note that we choose the form of $G$ using the matrix $J_m$ as in Section \ref{subsec_realforms}.

We define $\cG$ and the embedding of $G$ as follows:
\[
\begin{array}{ccl}
\cG      &:=& GL_{2m}(\bbC) \\
\cup     &  &                                    \\
G^\prime &:=& Sp_{2m}(\bbC) \cong G
\end{array}
\] with the standard embedding of $Sp_{2m}$.

Let $\cK_L$, $\cK_R$, $W_L$, $W_R$ be as in the $SO^*(2m)$ case. We
set $\cK := \cK_L \times \cK_R$, and embed in $\cG$ as:
\[
\begin{array}{rccl}
\cG     & \supset & \cK      & :=  \left\{ \left[ \begin{array}{cc} g & 0 \\ 0 & h \end{array}\right]: g \in \cK_L, h \in \cK_R \right\} \\
        &         & \cup     &                                     \\
K       & \cong   & K^\prime. &
\end{array}
\]
Let $\overline W$ be the standard representation of $\cG$, which we
also regard as the standard representation of $G$.  When restricted
to $\cK$, $\overline W \cong W_L \tot 1_\bbC \oplus 1_\bbC \tot
W_R$, and when restricted to $K$, $\overline W \cong W \oplus
W^*$, where $W$ is the standard representation of $K$.

As a representation of $G$, $\fg \cong \mcS^2(\cW)$. Let $\cmcS :=
\mcS(\cW)$ where $\mcS^2(\overline W)$, which we view as a graded,
representation of $\cG$. Using Theorem \ref{thm_SMn}, we have a
multiplicity free $\cG$-decomposition:
\[
    \mcS^d(\cW) \cong \bigoplus \F{2\la}{2m}
\] where the sum is over all partitions $\la$ such that $|\la| = d$ and $\ell(\la) \leq 2m$.

We now decompose relative to the action of $\cK$ using Proposition
\ref{prop_restrict_LRC}:
\[
    \mcS^d(\cW) \cong \bigoplus \c{2\la}{\mu}{\nu} \F{\mu}{m} \tot \F{\nu}{m}
\]
where $\mu$ and $\nu$ are partitions with $|\mu|+|\nu| = 2|\la|$ and
$\ell(\mu), \ell(\nu) \leq m$.  When we restrict to $K$,
\[
    \mcS^d(\cW) \cong \bigoplus \c{2\la}{\mu}{\nu} \F{\mu}{m} \ot \dual{\F{\nu}{m}}.
\]
We obtain a $K$ invariant exactly when $\mu = \nu$ as before. Equation \ref{eq_main_thm_Sp(m,R)}
follows.

\subsection{CASE: $G_0 = O(n, \bbC)$    }
\[
\begin{array}{|rcl|rcl|} \hline
G   &=&  O_n(\bbC) \times  O_n(\bbC) \hspace{15mm} & K   &\cong& O_n(\bbC) \\
\hline
\end{array}
\]
The group $K$ is diagonally embedded in $G$, as will be the case in
the $G_0 = Sp(m, \bbC)$ case.  We define $\cG$ and the embedding of
$G$ as follows:
\[
\begin{array}{ccl}
\cG      &:=& GL_n(\bbC) \times GL_n(\bbC)      \\
\cup     &  &                                         \\
G^\prime &:=& \{ (g,h)| g,h \in O_n(\bbC) \} \cong G
\end{array}
\] with the standard embedding of $O_n(\bbC)$ in $GL_n(\bbC)$.
Define: $\cK := \{ (g,g) | g \in GL_n(\bbC) \} \subset \cG$.

Set $\cW := \wedge^2(W) \tot 1_\bbC \oplus 1_\bbC \tot \wedge^2(W)$,
where $W$ is the standard representation of $GL_n(\bbC)$.  We view
$\cW$ as a representation of $\cG$, but as a representation of $G$,
$\fg \cong \cW$. Define:
\[\cmcS := \mcS(\cW) \cong \mcS(\wedge^2(W)) \tot
\mcS(\wedge^2(W)),
\]
which is a graded representation of $\cG$. Using Theorem
\ref{thm_SMn}, we have a multiplicity free $\cG$-decomposition:
\[
    \mcS^d(\cW) \cong \bigoplus \F{(2\mu)^\prime}{n} \tot
    \F{(2\nu)^\prime}{n}
\] where the sum is over all partitions $\mu$
and $\nu$ such that $|\mu| + |\nu| = d$ and $\ell((2\mu)^\prime)$,
$\ell((2\nu)^\prime)$ $\leq n$.

We now decompose relative to the action of $\cK$ using Proposition
\ref{prop_tensor_LRC}:
\[
    \mcS^d(\cW) \cong \bigoplus \c{\g}{(2\mu)^\prime}{(2\nu)^\prime} \F{\g}{n}
\] where $\mu$ and $\nu$ are as before and $\g$ is such that $|\g| =
2(|\mu|+|\nu|)$ and $\ell(\g)\leq n$.  This is the
$\cK$-decomposition.  By Theorem \ref{thm_Cartan-Helgason_GLnOn} we
see that we obtain a $K$-invariant exactly when $\g$ is of the form
$2 \la$ for some $\la$. Equation \ref{eq_main_thm_O(n,C)} follows.

\subsection{CASE: $G_0 = Sp(m, \bbC)$}
\[
\begin{array}{|rcl|rcl|} \hline
G   &=&  Sp_{2m}(\bbC) \times  Sp_{2m}{\bbC} \hspace{15mm} & K   &\cong& Sp_{2m}(\bbC) \\
\hline
\end{array}
\]
The group $K$ is diagonally embedded in $G$, as:
\[K \cong K^\prime := \{ (k,k) | k \in K \}.\]

We define $\cG$ and the embedding of $G$ as follows:
\[
\begin{array}{ccl}
\cG      &:=& GL_{2m}(\bbC) \times GL_{2m}(\bbC)      \\
\cup     &  &                                         \\
G^\prime &:=& \{ (g,h)| g,h \in Sp_{2m}(\bbC) \} \cong G
\end{array}
\] with the standard embedding of $Sp_{2m}(\bbC)$ in $GL_{2m}(\bbC)$.
Define: $\cK := \{ (g,g) | g \in GL_{2m}(\bbC) \} \subset \cG$.

Set $\cW := \mcS^2(\cW) \tot 1_\bbC \oplus 1_\bbC \tot \mcS^2(\cW)$,
where $W$ is the standard representation of $GL_{2m}(\bbC)$.  We
view $\cW$ as a representation of $\cG$, but as a representation of
$G$, $\fg \cong \cW$. Define:
\[\cmcS := \mcS(\cW) \cong \mcS(\mcS^2(W)) \tot
\mcS(\mcS^2(W)),
\]
which is a graded representation of $\cG$. Using Theorem
\ref{thm_SMn}, we have a multiplicity free $\cG$-decomposition:
\[
    \mcS^d(\cW) \cong \bigoplus \F{2\mu}{2m} \tot \F{2\nu}{2m}
\] where the sum is over all partitions $\mu$
and $\nu$ such that $|\mu| + |\nu| = d$ and $\ell(\mu), \ell(\nu)
\leq 2m$.

We now decompose relative to the action of $\cK$ using Theorem
\ref{prop_tensor_LRC}:
\[
    \mcS^d(\cW) \cong \bigoplus \c{\g}{2\mu}{2\nu} \F{\g}{2m}
\] where $\mu$ and $\nu$ are as before and $\g$ is such that $|\g| =
2(|\mu|+|\nu|)$ and $\ell(\g)\leq 2m$.  This is the
$\cK$-decomposition.  By Theorem \ref{thm_Cartan-Helgason_GL2mSp2m}
we see that we obtain a $K$-invariant exactly when $\g$ is of the
form $(2 \la)^\prime$ for some $\la$. Equation
\ref{eq_main_thm_Sp(m,C)} follows.

\subsection{CASE: $G_0 = O(p,q)$     }  Set $n := p+q$.
\[
\begin{array}{|rcl|rcl|} \hline
G   &=&   O_n(\bbC) \hspace{15mm} & K   &\cong& O_p(\bbC) \times O_q(\bbC)     \\
\hline
\end{array}
\]
The group $K$ is embedded in $G$ as a direct sum in the same way as
the $G_0 = U(p,q)$ case.

We define $\cG := GL_n(\bbC)$ and embed $G$ in the standard way. We
set
\[
\begin{array}{rcl}
\cK      &:=& \left\{ \left. \left[ \begin{array}{cc} g & 0 \\ 0 & h
\end{array} \right] \right| g \in GL_p(\bbC), h \in GL_q(\bbC)
\right\} \\
\cup & & \\
K^\prime &:=& \left\{ \left. \left[ \begin{array}{cc} g & 0 \\ 0 & h
\end{array} \right] \right| g \in O_p(\bbC), h \in O_q(\bbC)
\right\} \cong K,
\end{array}
\] where $K^\prime$ is defined to describe the embedding of $K$ in $\cK$.
Note that $O_p(\bbC)$ (resp. $O_q(\bbC)$) is embedded in $GL_p(\bbC)$ (resp. $GL_q(\bbC)$)
in the standard way.

Set $\cW := \wedge^2(W)$, where $W$ is the standard representation
of $\cG$.  As a representation of $G$, $\fg \cong \cW$. Define
$\cmcS := \mcS(\cW)$, which we view a graded representation of
$\cG$. Using Theorem \ref{thm_AMn}, we have a multiplicity free
$\cG$-decomposition:
\[
    \mcS^d(\cW) \cong \bigoplus \F{(2\la)^\prime}{n}
\] where the sum is over all partitions $\la$
such that $|\la| = d$ and $\ell((2\la)^\prime) \leq n$.

We now decompose relative to the action of $\cK$ using Theorem
\ref{prop_restrict_LRC}:
\[
    \mcS^d(\cW) \cong \bigoplus \c{(2\la)^\prime}{\a}{\b} \F{\a}{p} \tot \F{\b}{q}
\] where $\a$ and $\b$ are such that $d = 2|\la| =
|\a|+|\b|$ and $\ell(\a) \leq p$, $\ell(\b) \leq q$.  This is the
$\cK$-decomposition.  By Theorem \ref{thm_Cartan-Helgason_GLnOn} we
see that we obtain a $K$-invariant exactly when $\a$ (resp. $\b$) is
of the form $2 \mu$ (resp. $2 \nu$) for some $\mu$ (resp. $\nu$).
Equation \ref{eq_main_thm_O(p,q)} follows.

\subsection{CASE: $G_0 = Sp(p,q)$     }Set $m := p+q$, and $n = 2m$.
\[
\begin{array}{|rcl|rcl|} \hline
G   &=&   Sp_{2m}(\bbC) \hspace{15mm} & K   &\cong& Sp_{2p}(\bbC) \times Sp_{2q}(\bbC)     \\
\hline
\end{array}
\]
The group $K$ is embedded in $G$ as a direct sum in the same way as
the $G_0 = U_{p,q}$ and $G_0 = O(p,q)$ cases.  Note that we take
form defining $G$ to be the one given by the matrix $C_m$ as in Section \ref{subsec_realforms}.

We define $\cG := GL_{2m}(\bbC)$ and embed $G$ in the standard way. We
set
\[
\begin{array}{rcl}
\cK      &:=& \left\{ \left. \left[ \begin{array}{cc} g & 0 \\ 0 & h
\end{array} \right] \right| g \in GL_{2p}(\bbC), h \in GL_{2q}(\bbC)
\right\} \\
\cup & & \\
K^\prime &:=& \left\{ \left. \left[ \begin{array}{cc} g & 0 \\ 0 & h
\end{array} \right] \right| g \in Sp_{2p}(\bbC), h \in Sp_{2q}(\bbC)
\right\} \cong K,
\end{array}
\] where $K^\prime$ is defined to describe the embedding of $K$ in $\cK$, as in
the previous case.

Set $\cW := \mcS^2(W)$, where $W$ is the standard representation of
$\cG$.  As a representation of $G$, $\fg \cong \cW$. Define $\cmcS
:= \mcS(\cW)$, which we view a graded representation of $\cG$. Using
Theorem \ref{thm_SMn}, we have a multiplicity free
$\cG$-decomposition:
\[
    \mcS^d(\cW) \cong \bigoplus \F{2\la}{n}
\] where the sum is over all partitions $\la$
such that $|\la| = d$ and $\ell(\la) \leq n$.

We now decompose relative to the action of $\cK$ using Theorem
\ref{prop_restrict_LRC}:
\[
    \mcS^d(\cW) \cong \bigoplus \c{2\la}{\a}{\b} \F{\a}{2p} \tot \F{\b}{2q}
\] where $\a$ and $\b$ are such that $d = 2|\la| =
|\a|+|\b|$ and $\ell(\a) \leq 2p$, $\ell(\b) \leq 2q$.  This is the
$\cK$-decomposition.  By Theorem \ref{thm_Cartan-Helgason_GL2mSp2m}
we see that we obtain a $K$-invariant exactly when $\a$ (resp. $\b$)
is of the form $(2 \mu)^\prime$ (resp. $(2 \nu)^\prime$) for some
$\mu$ (resp. $\nu$). Equation \ref{eq_main_thm_Sp(p,q)} follows.

\section{Closed expressions}\label{sec_closed_expressions}
The purpose of this section is to show that in six of the ten cases,
we can actually close the sums involving Littlewood-Richardson
coefficients.  Since the ten cases naturally break into three
disjoint subsets, we are only required to address three of the sums.
Of these three cases, one has a simple and rather elegant
expression, another is considerably more complicated while it
remains an open problem to close the third.

This latter case, corresponds to the groups $O(n,\bbC)$,
$Sp(m,\bbC)$, $Sp(p,q)$ and $O(p,q)$.  This is a subject of further
study, but is not likely to close to a simple formula. The problem
would amount to closing:
\[
    \sum_{\la,\mu,\nu} \c{2\la}{(2\mu)^\prime}{(2\nu)^\prime} \, t^{|\la|}
\]
which, of course, is the same as:
\[
    \sum_{\la,\mu,\nu} \c{(2\la)^\prime}{2\mu}{2\nu} \, t^{|\la|}.
\]
Probably, obtaining a simple formula for either of the above
expressions would be almost as hard as:
\[
    \sum_{\la, \mu, \nu} \c{\la}{\mu}{\nu} \, t^{|\la|},
\]
which is a listed open problem on Richard Stanley's web page:
\begin{center}
\url{http://www-math.mit.edu/~rstan/ec/ch7supp.pdf}
\end{center}

Also on Richard Stanley's web page is:
\[
    \sum_{\la, \mu} \c{2\la}{\mu}{\mu} t^{|\la|} =
    \prod_{i \geq 1} \frac{1}{\sqrt{1-t^i}} \cdot
    \prod_{j \geq 1} \frac{1}{(1-t^{2j})^{2^{j-2}}}.
\]
Stanley's proof is this impressive formula is proved using some
rather technical manipulations of $q$-series, which we omit.  We do
however note:
\[
    \sum_{\la, \mu} \c{(2\la)^\prime}{\mu}{\mu} t^{|\la|} =
    \sum_{\la, \mu} \c{2\la}{\mu}{\mu} t^{|\la|}.
\]
From our point of view, the significance of the above expressions
are that they closes the sums of Littlewood-Richardson coefficients
that appear in the $GL(n,\bbR)$, $GL(m, \bbH)$, $SO^*(2n)$, and
$Sp(m,\bbR)$ cases.  For a combinatorial interpretation of the
coefficients see \cite{Will01}.

Now, as a final result, we will close the sum of
Littlewood-Richardson coefficients occurring in the $U(p,q)$ case,
which we have seen is the same as closing the sum in the
$GL(n,\bbC)$ case. The formula we obtain also appears on the above
web page, but we present here a representation theoretic proof as it
may be of independent interest.
\begin{thm}\label{thm_c^2-formula}
\[ \sum \left( \c{\la}{\mu}{\nu} \right)^2 t^{|\la|}
     = \prod_{k=1}^\infty {\frac{1}{1-2t^k}}.
\]
\end{thm}
In order to prove the above theorem, we will need to introduce some
additional notation.  Let $\vx = (x_1, x_2, \cdots)$, $\vy = (y_1,
y_2, \cdots)$, and $\vz = (z_1, z_2, \cdots)$ be three sets of
countably infinite indeterminates.  We now expand the following
symmetric product into Schur functions:
\begin{equation}\label{eq_g-def}
    \prod_{i,j,k=1}^\infty \left(\frac{1}{1-x_i y_j z_k} \right) =
    \sum_{\la, \mu, \nu} g_{\la \mu \nu }
    s_\mu(\vx) s_\nu(\vy) s_\la(\vz).
\end{equation}
The coefficients, $g_{\la \mu \nu}$, are non-negative integers which
can be interpreted as the tensor product multiplicities for the
symmetric groups (see \cite{Macd95, Stan99}).  That is to say, for
partitions $\la$, $\mu$ and $\nu$ of size $m$,
\[
    U_\mu \ot U_\nu \cong \bigoplus_{\la:|\la|=m} g_{\la \mu \nu} U_\la
\] where $U_\g$ is the irreducible $S_m$-representation indexed by
$\g$ (with $|\g|=m$) as in \cite{GW98}.  Note that we have $g_{\la
\mu \nu} = g_{\a \b \g}$ where $\a \b \g$ is a permutation of $\la
\mu \nu$.  This fact is a consequence of the fact that the
representations of $S_m$ are self dual.

By specialization $z_k = t^k$ in Equation \ref{eq_g-def} we obtain a
formal series in $\vx$ and $\vy$.  Define the coefficient of $s_\mu(\vx)
s_\nu(\vy)$ to be:
\[
    G_{\mu \nu }(t) := \sum_\la g_{\la \, \mu \, \nu} s_\la(t, t^2, t^3, \cdots).
\]
We will now provide an interpretation of $G_{\mu \nu }(t)$ involving
the representation theory of $GL_n(\bbC)$ as $n \rightarrow \infty$.
\begin{prop}\label{prop_from_stable_stanley}  For $\fg = \gl_n(\bbC)$ we
have:
\begin{equation}\label{eq_Stan84}
    \sum_{d=0}^\infty \left( \lim_{n \rightarrow \infty} \left[ \F{(\mu,\nu)}{n}, \mcS^d(\fg)
    \right] \right)
    t^d= \frac{G_{\mu \nu}(t)}{\prod_{k=1}^\infty \left( 1-t^k
    \right)}
\end{equation}
\end{prop}
\begin{proof}
    This result was first proved in \cite{Stan84}.  Please
    see Remark \ref{rem_Stan84} for an outline of the proof
    in the notation of this paper.
\end{proof}
The fact that the denominator of the right side of Equation
\ref{eq_Stan84} formally looks like the familiar $\eta$-function has
a clear interpretation. Given a connected, reductive, linear
algebraic group $G$ (over $\bbC$) with Lie algebra $\fg$, Kostant's
theorem asserts $\mcS(\fg) \cong \mcS(\fg)^G \otimes \mcH(\fg)$,
where $\mcH(\fg)$ denotes the space of harmonic polynomials on $\fg$
(see \cite{Kost63}).  We will see a reflection of this celebrated
theorem in the following discussion where we take $G = GL_n(\bbC)$
and $\fg = \gl_n(\bbC)$.

Set $n = p+q$ and embed $K = GL_p(\bbC) \times GL_q(\bbC)$ in $G =
GL_n(\bbC)$ as in the proof of Equation \ref{eq_main_thm_U(p,q)}. We
will now identify the distribution of the $K$-invariants in the
harmonic polynomials on $\fg = \gl_n(\bbC)$. Define:
\[
h_{p,q}(d):=\dim\left(\mcH^d(\fg)\right)^{GL_p \times GL_q}
\]
and
\[
f_{p,q}(d):=\dim\left(\mcS^d(\fg)\right)^{GL_p \times GL_q}.
\]
From Equation \ref{eq_main_thm_U(p,q)} we see that $f_{p,q}(d) =
f_{p,q}( d_0)$ for all $d \geq d_0$ where $d_0 = \min(p,q)$.  As we
remarked before, this fact justifies the definition
\[
f(d) := \lim_{p,q \rightarrow \infty} f_{p,q}(d).
\]
We also set $F(t) = \sum_{d=0}^\infty f(d) t^d$. Proposition
\ref{eq_main_thm_U(p,q)} implies that:
\[ F(t) = \sum \left( \c{\la}{\mu}{\nu} \right)^2 t^{|\la|} \]
where the sum is over all triples of non-negative integer partitions
$\la$, $\mu$ and $\nu$.

Let $I_n(t) = \sum_{d=0}^\infty \dim \mcS^d(\fg)^{GL_n(\bbC)} \;
t^d$. It is a standard fact due to Chevalley (see \cite{Chev55})
that $\mcS(\fg)^G$ is a polynomial ring.  For the case of
$GL_n(\bbC)$, we can choose generators for this ring in degrees
$1,2, \cdots, n$. Thus $I_n(t) = \prod_{k=1}^n \frac{1}{1-t^k}$.
Note that as $n \rightarrow \infty$ the coefficients of the product
stabilize. In light of this fact, we set $ I(t) :=
\prod_{k=1}^\infty \left( \frac{1}{1-t^k} \right)$.

A combination of this fact and Kostant's theorem imply that $F(t) =
I(t) H(t)$, where $H(t)$ is a formal series in $t$ with non-negative integer
coefficients. Theorem
\ref{thm_c^2-formula} follows from:
\begin{prop}
\[    H(t) = \prod_{k=0}^\infty \left( \frac{1-t^k}{1-2 t^k}
\right) \]
\end{prop}
\begin{proof}
For $d\geq 0$, define $h(d)$ by $H(t) = \sum_{d=0}^\infty h(d) t^d$.
From the stability properties of $I_n(t)$ and $f_{p,q}(d)$ it is
easy to see that $h(d) = \lim_{p,q \rightarrow \infty} h_{p,q}(d)$.

Applying the Cartan-Helgason Theorem from Section
\ref{thm_Cartan-Helgason_GLnGLpGLq} and Proposition
\ref{prop_from_stable_stanley} we obtain $H(t) = \sum_\mu G_{\mu
\mu}(t)$. Therefore, we will consider the symmetric function
\[\sum_{\la, \mu} g_{\la \, \mu \, \mu} s_\la(\vz)\]
and evaluate at the point $\vz = (t, t^2, t^3, \cdots)$ thus
obtaining $H(t)$.  We now interpret the above sum in term of the
representation theory of the symmetric group.

As a representation of $S_m \times S_m$ under left and right
multiplication we have:
\[ \bbC[S_m] \cong \bigoplus_{\mu \vdash m} U_\mu \ot U_\mu, \]
(where $\bbC[S_m]$ is the group algebra of the symmetric group). If
we restrict this action to the diagonal subgroup we obtain: $
\bbC[S_m] \cong \bigoplus_{\la \vdash m} \left( \sum_\mu g_{\la \,
\mu \, \mu} \right) U_\la. $ Observe that the diagonal action of
$S_m$ on $\bigoplus  U_\mu \ot U_\mu$ is, as a representation,
isomorphic to the conjugation action of $S_m$ on $\bbC[S_m]$.

Given a representation $\rho: S_m \rightarrow GL(V)$, denote the
Frobenius characteristic (see \cite{Macd95, Stan99}) of $V$ by
$\ch(V)$. That is,
\[
    \ch(V)(\sigma) = \frac{1}{m!} \sum_{\sigma \in S_m}
    \chi_\rho(\sigma) p_{\la(\sigma)}(\vx),
\]
where $p_\mu(\vx)$ denotes the power symmetric function, which for a
partition $\mu$, is defined as $p_\mu(\vx) =
p_{\mu_1}(\vx)p_{\mu_2}(\vx)\cdots$ with $p_k(\vx) = x_1^k + x_2^k +
\cdots$, and $\la(\sigma)$ is a partition denoting the shape of
$\sigma$ in its disjoint cycle notation  (see \cite{Macd95}).

Suppose $V = \bbC[S_m]$ with $S_m$ acting by conjugation. For a
permutation $\sigma$, let $C_{S_m}(\sigma)$ denote the centralizer
of $\sigma$ in $S_m$. Then,
\begin{eqnarray*}
    \ch ( \bbC[S_m] ) &=& \frac{1}{m!} \sum_{\sigma \in S_m} |C_{S_m}(\sigma)| \; p_{\la(\sigma)}(\vx)
    = \sum_{\mu \vdash m} z_\mu \; \frac{p_\mu(\vx)}{z_\mu} = \sum_{\mu \vdash m}
    p_\mu(\vx).
\end{eqnarray*} where $z_\mu = 1^{a_1} 2^{a_2} \cdots a_1! a_2!
\cdots$ with $\mu$ having $a_i$ cycles of size $i$. \\
\noindent Summing these symmetric functions together we obtain:
\begin{equation*}
    \sum_{m=0}^{\infty} \ch( \bbC[S_m] ) = \sum_\mu p_\mu(\vx)
                                          = \prod_{k=1}^\infty
                                          \frac{1}{1-p_k(\vx)}.
\end{equation*}
On the other hand, note that $\ch(U_\la) = s_\la(\vx)$ and therefore
we have shown that:
\begin{equation*}
    \sum_{\la, \mu} g_{\la \, \mu \, \mu} s_\la(\vx) = \prod_{k=1}^\infty \frac{1}{1-p_k(\vx)}.
\end{equation*}
And so,
\begin{eqnarray*}
    H(t) &=& \left(\prod_{k=1}^\infty \frac{1}{1-p_k}\right)(t, t^2, t^3,
    \cdots)= \prod_{k=1}^\infty \frac{1}{1-p_k(t, t^2, t^3, \cdots)} = \prod_{k=1}^\infty \frac{1}{1-(t^k+t^{2k}+ \cdots)} \\
    &=& \prod_{k=1}^\infty \frac{1}{1-\frac{t^k}{1-t^k}} = \prod_{k=1}^\infty \frac{1-t^k}{1-2t^k}.
\end{eqnarray*}
\end{proof}
\begin{rem}
Note that the coefficient of $t^n$ in $\prod_{k=1}^\infty
\frac{1-t^k}{1-q t^k}$ is the number of conjugacy classes in the
finite group $GL(n, F)$ where $F$ is the field with $q$ elements
(see \cite{BHF86}).
\end{rem}
\begin{rem}\label{rem_Stan84}
We briefly describe the proof of Proposition
\ref{prop_from_stable_stanley} from a point of view closer to the
representation theory of classical groups as described in
\cite{HTW05}. First, we note from Corollary \ref{cor_dual_Mpq} that
we have:
\[
    \mcS^d(\fg) \cong \bigoplus_{\twoline{\rho: |\rho| = d}{\ell(\rho)\leq n}}
    \dual{\F{\rho}{n}} \ot \F{\rho}{n}
\] as a representation of $GL_n(\bbC)$.  In our notation, $\dual{\F{\rho}{n}} \cong
\F{(0,\rho)}{n}$.

Using the Theorem 2.1.1 of \cite{HTW05} we obtain:
\begin{equation}\label{eq_HTW05}
    \lim_{n \rightarrow \infty}
    \left[ \F{(0,\rho)}{n} \ot \F{(\rho,0)}{n}, \F{(\mu,\nu)}{n}
    \right] = \sum \c{\rho}{\la}{\mu} \c{\rho}{\la}{\nu}
\end{equation}
where the sum is over all partitions $\la$.  (Note that the sum is
easily seen to be finite after taking into account the support of
the Littlewood-Richardson coefficients.)  This formula perhaps first
appeared \cite{King71} as (4.6) with (4.15).

The partitions $\rho$, $\mu$ and $\nu$ are fixed on the right side
of Equation \ref{eq_HTW05}.  We will put together a formal sum of
Schur functions to encode these numbers combinatorially as:
\[
    S(\vx,\vy; t) = \sum_{\mu,\nu,\rho} \left( \sum_\la \c{\rho}{\la}{\mu} \c{\rho}{\la}{\nu}
    \right) s_\mu(\vx) s_\nu(\vy) \, t^{|\rho|} = \sum_{\rho,\la} s_{\rho/\la}(\vx)
    s_{\rho/\la}(\vy) \, t^{|\rho|}
\] where we define the ``Skew Schur" function as in \cite{Macd95, Stan99}:
\[
    s_{\a/\b}(\vx) = \sum_\g \c{\a}{\b}{\g} s_\g(\vx).
\]
Proposition \ref{prop_from_stable_stanley} now follows from the
identity:
\[
    \sum_{\rho,\la} s_{\rho/\la}(\vx) \, s_{\rho/\la}(\vy)
    t^{|\rho|} = \prod_k \left[\left(\frac{1}{1-t^k}\right) \prod_{i,j} \left(\frac{1}{1 - x_i
    y_j}\right)\right],
\]  which can be found in \cite{Macd95} p. 94 (28a).

Lastly, we point out that several analogous skew identities in
\cite{Macd95, Stan99} may be interpreted in the context of the
generalization of Kostant's result in \cite{KR71} as in
\cite{Will02}.
\end{rem}

\def\cprime{$'$} \def\cprime{$'$}


\begin{bibdiv}

\begin{biblist}

\bib{BHF86}{article}{
    author={Benson, David},
    author={Feit, Walter},
    author={Howe, Roger},
     title={Finite linear groups, the Commodore 64, Euler and Sylvester},
   journal={Amer. Math. Monthly},
    volume={93},
      date={1986},
    number={9},
     pages={717\ndash 719},
      issn={0002-9890},
    review={MR863974 (87m:05011)},
}

\bib{Chev55}{article}{
    author={Chevalley, C.},
     title={Sur certains groupes simples},
   journal={T\^ohoku Math.J.(2)},
    volume={7},
      date={1955},
     pages={14\ndash 66},
      issn={0040-8735},
    review={MR0073602 (17,457c)},
}

\bib{GW98}{book}{
      author={Goodman, R.},
      author={Wallach, N.R.},
       title={Representations and invariants of the classical groups},
   publisher={Cambridge University Press},
     address={Cambridge},
        date={1998},
        ISBN={0-521-58273-3},
      review={\MR{99b:20073}},
}

\bib{Howe89-1}{article}{
    author={Howe, Roger},
     title={Remarks on classical invariant theory},
   journal={Trans. Amer. Math. Soc.},
    volume={313},
      date={1989},
    number={2},
     pages={539\ndash 570},
      issn={0002-9947},
    review={MR986027 (90h:22015a)},
}

\bib{Howe89-2}{article}{
    author={Howe, Roger},
     title={Transcending classical invariant theory},
   journal={J. Amer. Math. Soc.},
    volume={2},
      date={1989},
    number={3},
     pages={535\ndash 552},
      issn={0894-0347},
    review={MR985172 (90k:22016)},
}

\bib{Howe95}{incollection}{
      author={Howe, Roger},
       title={Perspectives on invariant theory: {S}chur duality,
  multiplicity-free actions and beyond},
        date={1995},
   booktitle={The schur lectures (1992) (tel aviv)},
   publisher={Bar-Ilan Univ.},
     address={Ramat Gan},
       pages={1\ndash 182},
      review={\MR{96e:13006}},
}

\bib{HTW05}{article}{
    author={Howe, Roger},
    author={Tan, Eng-Chye},
    author={Willenbring, Jeb F.},
     title={Stable branching rules for classical symmetric pairs},
   journal={Trans. Amer. Math. Soc.},
    volume={357},
      date={2005},
    number={4},
     pages={1601\ndash 1626 (electronic)},
      issn={0002-9947},
    review={MR2115378 (2005j:22007)},
}

\bib{King71}{article}{
    author={King, R. C.},
     title={Modification rules and products of irreducible representations
            of the unitary, orthogonal, and symplectic groups},
   journal={J. Mathematical Phys.},
    volume={12},
      date={1971},
     pages={1588\ndash 1598},
    review={MR0287816 (44 \#5019)},
}


\bib{Kost63}{article}{
    author={Kostant, B.},
     title={Lie group representations on polynomial rings},
   journal={Amer. J. Math.},
    volume={85},
      date={1963},
     pages={327\ndash 404},
      issn={0002-9327},
    review={MR0158024 (28 \#1252)}
}

\bib{KR71}{article}{
      author={Kostant, B.},
      author={Rallis, S.},
       title={Orbits and representations associated with symmetric spaces},
        date={1971},
     journal={Amer. J. Math.},
      volume={93},
       pages={753\ndash 809}
}

\bib{Macd95}{book}{
      author={Macdonald, I.G.},
       title={Symmetric functions and {H}all polynomials},
     edition={Second},
   publisher={The Clarendon Press Oxford University Press},
     address={New York},
        date={1995},
        ISBN={0-19-853489-2},
        note={With contributions by A. Zelevinsky, Oxford Science
  Publications},
      review={\MR{96h:05207}},
}

\bib{Stan84}{article}{
      author={Stanley, Richard~P.},
       title={The stable behavior of some characters of {${\rm SL}(n,{\bf
  C})$}},
        date={1984},
        ISSN={0308-1087},
     journal={Linear and Multilinear Algebra},
      volume={16},
      number={1-4},
       pages={3\ndash 27},
      review={\MR{86e:22025}},
}

\bib{Stan97}{book}{
      author={Stanley, Richard~P.},
       title={Enumerative combinatorics. {V}ol. 1},
   publisher={Cambridge University Press},
     address={Cambridge},
        date={1997},
        ISBN={0-521-55309-1},
      review={\MR{98a:05001}},
}

\bib{Stan99}{book}{
      author={Stanley, Richard~P.},
       title={Enumerative combinatorics. {V}ol. 2},
   publisher={Cambridge University Press},
     address={Cambridge},
        date={1999},
        ISBN={0-521-56069-1},
      review={\MR{1 676 282}},
}

\bib{Will00}{thesis}{
      author={Willenbring, Jeb~F.},
       title={Stability properties for $q$-multiplicities and branching
  formulas for representations of the classical groups},
        type={Ph.D. Thesis},
        date={2000},
}

\bib{Will01}{article}{
    author={Willenbring, Jeb~F.},
     title={A stable range for dimensions of homogeneous ${\rm
            O}(n)$-invariant polynomials on the $n\times n$ matrices},
   journal={J. Algebra},
    volume={242},
      date={2001},
    number={2},
     pages={691\ndash 708},
      issn={0021-8693},
    review={MR1848965 (2002f:13013)},
}

\bib{Will02}{article}{
    author={Willenbring, Jeb~F.},
     title={An application of the Littlewood restriction formula to the
            Kostant-Rallis theorem},
   journal={Trans. Amer. Math. Soc.},
    volume={354},
      date={2002},
    number={11},
     pages={4393\ndash 4419 (electronic)},
      issn={0002-9947},
    review={MR1926881 (2003e:20049)},
}

\end{biblist}
\end{bibdiv}
\end{document}